\documentclass[10pt]{article}
\usepackage[T1]{fontenc} 
\usepackage[utf8]{inputenc}

\usepackage{lmodern}

\usepackage{calrsfs}
\usepackage{amssymb}
\usepackage[cmex10]{amsmath}
\usepackage{amsthm}
\usepackage{latexsym}
\usepackage{epsfig}
\usepackage{graphicx}

\DeclareMathOperator*{\maxobj}{\textup{maximize\,}}
\DeclareMathOperator*{\minobj}{\textup{minimize\,}}
\newcommand{\Prob}{\ensuremath{\mathsf{P}}}
\newcommand{\Ex}{\ensuremath{\mathsf{E}}}

\newtheorem{theorem}{Theorem}[section]

\newtheorem{lemma}[theorem]{Lemma}
\newtheorem{conjecture}[theorem]{Conjecture}
\newtheorem{proposition}[theorem]{Proposition}

\usepackage{appendix}

\usepackage[top=1in, bottom=1.25in, left=0.9in, right=0.9in]{geometry}

\usepackage[sort&compress]{natbib} 
\setcitestyle{authoryear,open={(},close={)}}

\usepackage[colorlinks=true,breaklinks=true,bookmarks=true,urlcolor=blue,
     citecolor=blue,linkcolor=blue,bookmarksopen=false,draft=false]{hyperref}

\def\EMAIL#1{\href{mailto:#1}{#1}}
\def\URL#1{\href{#1}{#1}}         

\begin{document}

\title{Towards minimum loss job routing to parallel heterogeneous multiserver queues via index policies}


\author{Jos\'e Ni\~no-Mora 
\\ Department of Statistics \\
    Carlos III University of Madrid \\
     28903 Getafe (Madrid), Spain \\  \EMAIL{jnimora@alum.mit.edu}, \URL{http://orcid.org/0000-0002-2172-3983} }
\date{Published in \textit{European Journal of Operational Research}, vol.\ 220, 705--715, 2012 \\
DOI: \href{https://doi.org/10.1016/j.ejor.2012.02.033}{10.1016/j.ejor.2012.02.033}
}

\maketitle

\begin{abstract}%
This paper considers a Markovian model for the optimal dynamic routing of homogeneous traffic to parallel heterogeneous queues, each having its own finite input buffer and server pool, where buffer and server-pool sizes, as well as service rates, may differ across queues. The main goal is to identify a heuristic index-based routing policy with low complexity that consistently attains a nearly minimum average loss rate (or, equivalently, maximum throughput rate). A second goal is to compare alternative policies, with respect to computational demands and empirical performance. A novel routing policy that can be efficiently computed is developed based on a second-order extension to Whittle’s restless bandit (RB) index, since the latter is constant for this model. New results are also given for the more computationally demanding index policy obtained via policy improvement (PI), including that it reduces to shortest queue routing under symmetric buffer and server-pool sizes. A numerical study shows that the proposed RB index policy is nearly optimal across the instances considered, and substantially outperforms several previously proposed index policies.\end{abstract}%

\newpage
\section{Introduction}
\label{s:intro}
This paper considers a Markov decision process (MDP) model for a distributed service system, concerning the optimal  dynamic routing, or assignment, of a Poisson job arrival stream with rate $\lambda$ to $K$ parallel multiserver queues, where the latter have finite input buffer capacities $n_k$ for queue $k \in \mathbb{K} \triangleq \{1, \ldots, K\}$.
Arrivals to a nonfull system are immediately and  irrevocably routed to a queue, while arrivals to a full system are blocked and lost. The performance objective is to minimize the average loss rate, i.e., the average number of arriving jobs per unit time that are blocked and lost.  
Each queue $k$ has its own pool of $m_k$ identical exponential servers, each working at rate $\mu_k$. 
The system is heterogeneous, in that buffer and server-pool sizes, as well as service rates, may differ across queues.

Note that the performance objective of minimizing the average loss rate is equivalent to that of maximizing the average throughput rate, i.e., the average number of jobs completed per unit time, since under any job assignment policy the sum of the throughput rate and the loss rate equals the arrival rate $\lambda$ (since a job that is not lost eventually completes service).

Such a model and performance objective are relevant to a variety of applications: 
(1) The main motivating application is \emph{load sharing} or \emph{balancing} to distribute client requests across a  cluster of Web-server nodes (see, e.g.,  \citet{ZhangFan08}). As for the motivation for considering node queues with finite buffers, the latter yield a simple form of \emph{admission control}  (as jobs finding all buffers full are rejected), which the service provider can leverage to offer statistical response-time  Quality of Service (QoS) guarantees  in a Service Level Agreement for admitted requests. 
The main performance metrics in a cluster are the average response time  and the throughput.
Since in such systems the response times are controlled by limiting node buffer sizes, which has the undesirable side effect of causing request losses, 
the relevant optimization problem is to minimize the average  loss rate.
This is in contrast to most analytical work on Web-server load sharing, which mainly addresses response-time minimization, often assuming infinite buffer space at the nodes.
(2) In  geographically distributed telecommunication networks, messages are relayed to one of a number of service centers. Buffers sizes at the latter are limited to ensure target 
QoS levels on response-time delays. 
See, e.g., \cite{PourSeidm92}.
(3) In high-speed packet-switched networks,  packets flowing through
  a router are forwarded through one of multiple outgoing finite-buffer channels. When all buffers are full, packets are dropped and lost. As noted in \citet{kangTan97}, ``packet loss due to buffer overflow is the primary performance measure in high-speed networks \ldots''.

Obtaining an optimal policy for the model of concern here by numerical solution of its dynamic programming (DP) equations is, generally, a computationally intractable task for realistic model sizes. 
Thus, e.g., in a cluster with $K$ nodes where each node has $n$ buffer spaces, the DP equations have $(n+1)^K$ variables, one per system state. A well-known approach to tackle such  equations is to solve a related \emph{linear programming} (LP) problem, which has about $K (n+1)^K$ constraints on $(n+1)^K$ variables.
Even with a buffer size of $n = 60$ jobs at each node, and considering a very small cluster with $K = 3$ nodes, one gets an LP problem with nearly a quarter million variables and three quarter million constraints, which would stress current high-end workstations. Besides, the optimal policy would typically be too complex to implement.  Currently, sizes of dozens of servers are common, with many firms operating large-scale clusters 
 with hundreds or thousands of servers. The resulting DP equations would hence be impossible to solve numerically.

This motivates our main goal, which is to identify tractable heuristic policies that perform well. 
Particularly appealing is the class of \emph{index policies}. These attach a \emph{load index} $\theta_k(x_k)$ to each queue $k$, which is a numeric function of its \emph{load state} $x_k$, where the latter gives the number of jobs  in the queue input buffer, waiting or in service. The policy routes each arrival to a nonfull queue, if any, with currently lowest index value.
Classic examples are the \emph{Shortest (nonfull) Queue} (SQ) and the \emph{Shortest (nonfull) Expected Delay} (SED) routing rules.
Index policies,  being based on dynamic state information, offer the potential to significantly improve system performance over static (state-blind) policies,  while remaining tractable  if the indices can be evaluated with low complexity.

\subsection{Prior Work on Minimum Loss Routing to Finite-Buffer Queues}
\label{s:rl}
While a vast literature addresses problems of job assignment to parallel infinite-buffer queues to minimize the average response time, there is a gap in the literature on models for job assignment to parallel finite-buffer queues to minimize the average loss rate (maximize the average throughput rate), which have received much less attention. 
Early research established the optimality of simple index rules  under strong assumptions. 
For the case $m_k \equiv n_k \equiv 1$ of single-server queues without waiting room, the \emph{Fastest Available Server} (FAS)  rule, with index $\theta_k^{\textup{FAS}}(x_k) \triangleq 1/\mu_k$, is shown to be optimal in
\citet{seth77} and in
\citet{dermanetal80}.
\citet{hoko90} shows that SQ routing is optimal in the case of parallel single-server queues with equal service rates, allowing unequal buffers. 

More recently, the research focus has shifted to the design of suboptimal heuristic  policies with low complexity for intractable models with heterogeneous queues, such as the present model, for which, to the best of the author's knowledge, no analytical results are available on the structure of optimal policies.  
A well-known approach to obtain a static policy is the \emph{optimal Bernoulli Splitting} (OBS) method, which 
splits the arriving job stream into independent substreams assigned to the queues according to fixed probabilities, obtained by solving a \emph{nonlinear program} (NLP). 
Note that  such policies  may route a job to a  full queue, thus losing it, even when some other queue has free space.
\citet{yaoshant87} addresses computation of the OBS for maximum throughput routing to parallel multiserver  queues without waiting room ($m_k \equiv n_k$).

As for the design of more sophisticated index policies, a widely applied method is to start with the OBS, and then to carry out one step of the 
\emph{policy improvement} (PI) algorithm for MDPs, which is known to yield a better policy. 
Such a \emph{one-step PI method} is proposed in \citet{krish88} to obtain a load index $\theta_k^{\textup{PI}}(x_k)$ for the present model, albeit bypassing the issue of computing the OBS by referring to \citet{yaoshant87}, which only considers the no-waiting-room case.

Another powerful method to design index policies is to deploy  \emph{restless bandit} (RB) \emph{indexation}. See \citet{whit88}.
\citet{nmmp02} first applied  the
\emph{Whittle RB index} $\theta_k^{\textup{RB}}(x_k)$ to obtain policies for admission control and routing  to parallel queueing systems,   presenting also index extensions, sufficient conditions for \emph{indexability} (existence of the index), and an efficient index algorithm. 
The application of such results to various routing models  is outlined in  \citet{nmnetcoop07}.

\subsection{Goals and Contributions}
\label{s:contr}
The main goal of this paper is to identify a routing policy of index type for the model of concern that has low computational complexity and consistently attains a nearly minimum average loss rate.
A related second goal is to compare several alternative routing policies, elucidating their relative standing with respect both to computational demands and average loss performance.
 
The main contributions are the following:
\begin{enumerate}
\item Design of a new routing policy based on an extension to RB indexation. Whereas the Whittle index turns out to be noninformative for the present model, taking  on a constant value, a novel tie-breaking second-order routing index is obtained based on analysis of a related discounted problem. Such an index is efficiently computed with linear complexity on the buffer size by a first-order linear recursion, which is solved in closed form. Also, a closed formula is presented for the upper bound on the maximum throughput performance, based on an RB relaxation.
A benchmarking study shows that the proposed RB index policy is nearly optimal in every instance considered, while four alternative previously proposed index policies  are found to be severely suboptimal in moderate traffic.
\item New results on the OBS policy, which extend work in \citet{yaoshant87}, as the latter applies only to the no-waiting-room case. Algorithmic  results are presented for the efficient numerical solution of the NLP yielding the OBS. 
It is further shown that, under symmetric  buffer and server-pool sizes, the OBS
feeds all queues with equal offered loads.
\item New results on the PI index policy introduced in \citet{krish88}.
It is shown that, in the case of symmetric  buffer and server-pool sizes  across queues,  the PI index policy reduces to SQ routing, thus failing to take into account differences in service rates across queues. 
\end{enumerate}

\subsection{Organization of the paper}
\label{s:op}
The remainder of the paper is organized as follows.
After the model is described in Section \ref{s:mod}, Section \ref{s:rbim} develops the novel second-order RB routing index. 
Section \ref{s:obspim} discusses application of the PI method, extending previous work as outlined above.
Section \ref{s:bs} reports on the benchmarking study.
Section \ref{s:concl} concludes.

\section{Model Description}
\label{s:mod}
Consider a distributed service system consisting of
 $K$ parallel service stations, which are used to fulfill the service requests, or jobs,  submitted by a homogeneous customer stream, 
whose arrival times form a Poisson process with rate $\lambda$. 
Station $k \in \mathbb{K}$ has $m_k \geqslant 1$  parallel exponential servers, each working at rate 
$\mu_k$, and a  finite input buffer with capacity for holding $n_k$ jobs, waiting or being served, with $n_k \geqslant m_k$.
We will refer to station $k$ as \emph{queue} $k$ in the sequel.
Upon arrival of a job, a central controller irrevocably routes the job to a nonfull queue, if any is available, according to a \emph{routing policy}. 
If all buffers are full, an arriving job is blocked and lost.
Once in a queue, jobs are scheduled for service on a First-Come First-Served (FCFS) basis. 

Denote by $X_k(t) \in \mathbb{X}_k \triangleq \{0, 1, \ldots, n_k\}$ the \emph{state}
 of queue $k$ at time $t$,  giving the number of jobs it holds either waiting or being served, and by $A_k(t) \in \{0,
1\}$ the binary \emph{action}  that \emph{equals $0$  if
an arrival at time $t$  would be routed to queue}
$k$. It is convenient to imagine that each queue $k$ has its own \emph{entry gate}, which the system controller can open ($A_k(t) = 0$) or shut ($A_k(t) = 1$) at each time.
The requirement that no job should be routed to a full queue is thus formulated as
\begin{equation}
\label{eq:ak1}
A_k(t) = 1 \, \text{ if } \,  X_k(t) = n_k, \quad k \in \mathbb{K}, \quad t \geqslant 0,
\end{equation}
while the requirement that a job must be routed to a nonfull queue if one is available is  formulated as
\begin{equation}
 \label{eq:spconstr}
 \sum_{k \in \mathbb{K}} A_k(t) = K-1  \, \, \text{if} \, \, \mathbf{X}(t) \neq \mathbf{n}, \quad t \geqslant 0,
 \end{equation}
where $\mathbf{X}(t) = \big(X_k(t)\big)$ and $\mathbf{n} = (n_k)$.
Further, it is required that the number of jobs in each queue does not exceed its buffer capacity: 
\begin{equation}
\label{eq:bc}
X_k(t) \leqslant n_k, \quad k \in \mathbb{K}, \quad t \geqslant 0.
\end{equation}

Actions are selected according to a \emph{routing policy} $\boldsymbol{\pi}$ from the class $\boldsymbol{\Pi}$ of  
\emph{history-dependent randomized policies} (in standard MDP terminology), which choose the joint action $\mathbf{A}(t) = \big(A_k(t)\big)$ at each time $t \geqslant 0$ based on the history of joint states $\{\mathbf{X}(s)\colon 0 \leqslant s \leqslant t\}$ and actions $\{\mathbf{A}(s)\colon 0 \leqslant s < t\}$,  and  satisfy  the sample-path constraints (\ref{eq:ak1})--(\ref{eq:bc}).
We will further refer to the class $\boldsymbol{\Pi}^{\textup{SR}} $of \emph{stationary randomized} policies, where the joint action $\mathbf{A}(t)$ selected at each time $t$ depends only on the joint state $\mathbf{X}(t)$, possibly in a randomized fashion.
Under a \emph{stationary deterministic} policy, $\mathbf{A}(t)$ is a deterministic function of $\mathbf{X}(t)$.

Consider the 
\emph{optimal routing problem} of
finding an \emph{average-reward optimal}  policy that maximizes the average throughput rate.
The formulation of such a problem will be simplified by the fact that the corresponding MDP model has finite state and action spaces, and, besides, is  \emph{unichain}, meaning that the Markov chain generated on the state space by any stationary policy has a single ergodic class. 
This follows from the intuitive observation that, in the present model, every pair of system states communicates with positive probability under any stationary policy that satisfies (\ref{eq:ak1})--(\ref{eq:bc}).
The unichain property allows us to focus, without loss of generality, on routing policies from the class $\boldsymbol{\Pi}^{\textup{SR}}$. For any such policy, we will denote by 
$\widetilde{\mathbf{X}} = (\widetilde{X}_k)$ and $ \widetilde{\mathbf{A}} = (\widetilde{A}_k)$ random vectors with 
the steady-state distributions of the ergodic joint state process $\mathbf{X}(t)$ and joint action process $\mathbf{A}(t)$.
We can hence formulate the throughput maximization problem as
\begin{equation}
\label{eq:etmp2}
\maxobj_{\boldsymbol{\pi} \in \boldsymbol{\Pi}^{\textup{SR}}} 
\Ex^{\boldsymbol{\pi}}\left[\sum_{k \in \mathbb{K}}  \bar{\mu}_k\big(\widetilde{X}_k\big)\right],
\end{equation}
where 
$\Ex^{\boldsymbol{\pi}}[\cdot]$ denotes expectation under
policy $\boldsymbol{\pi}$, 
since the unichain property ensures that the maximum throughput rate  is independent of the initial state.

We will draw in the sequel on the observation that problem (\ref{eq:etmp2}) is equivalent to the problem of minimizing the average loss rate in the system, which, by the \emph{PASTA} (Poisson Arrivals See Time Averages) property, can be formulated as
\begin{equation}
\label{eq:pa2a}
\minobj_{\boldsymbol{\pi} \in \boldsymbol{\Pi}^{\textup{SR}}} 
\lambda \, \Prob^{\boldsymbol{\pi}}\left\{\widetilde{\mathbf{X}} = \mathbf{n}\right\},
\end{equation}
where $\Prob^{\boldsymbol{\pi}}\left\{\cdot\right\}$ denotes probability under policy $\boldsymbol{\pi}$.

The reason for the equivalence of problems (\ref{eq:etmp2}) and (\ref{eq:pa2a}) is that, under any routing policy $\boldsymbol{\pi} \in \boldsymbol{\Pi}^{\textup{SR}}$,  one has the invariance identity
\begin{equation}
\label{eq:idthloss}
\Ex^{\boldsymbol{\pi}}\left[\sum_{k \in \mathbb{K}}  \bar{\mu}_k\big(\widetilde{X}_k\big)\right] + \lambda \Prob^{\boldsymbol{\pi}}\left\{\widetilde{\mathbf{X}} = \mathbf{n}\right\} \equiv \lambda,
\end{equation}
i.e., the throughput rate plus the loss rate equals the arrival rate.

We will denote by $J^*$ the minimum loss rate objective for problem (\ref{eq:pa2a}), and by $V^*$ the maximum throughput rate objective for problem (\ref{eq:etmp2}). Note that $J^* + V^* = \lambda$.

 \section{RB Index Method}
\label{s:rbim}
The conventional RB index method is based on formulating (\ref{eq:etmp2}) as a \emph{multiarmed restless bandit problem} (MARBP), and then 
using the  index policy proposed in \citet{whit88} for the latter, which further yields a performance bound.
Here, we need to consider a variant of the MARBP to formulate the optimal routing problem of concern.
While the RB approach was first applied to the design of routing policies in \citet[Sec.\ 8.1]{nmmp02},  it is shown below that direct application of the results there to the present model does not yield a
specific policy,  since the Whittle index takes on a constant value.
As a tie-breaker, a   second-order  index is presented, based on analysis of a related discounted problem.

\subsection{Formulation as a Variant of the MARBP}
\label{s:fmarbp}
The continuous-time average-reward MARBP concerns the optimal dynamic allocation of effort to $K$ 
stochastic projects with state processes $X_k(t)$, for $k \in \mathbb{K}$, modeled as independent RBs, i.e., binary-action ($A_k(t) = 1$: active; $A_k(t) = 0$: passive)
MDPs. 
\citet{whit88} assumes that a fixed number $M$ out of the $K$ 
projects must be active at each time.
Separable rewards are earned at the state- and action-dependent rate $R_k(x_k, a_k)$ for project $k$. 
Denoting by $\boldsymbol{\Pi}^{\textup{SR}}(M)$ the  class of admissible \emph{scheduling policies} that are stationary randomized, satisfy the sample-path activity constraint $\sum_k A_k(t) \equiv M$, and generate an ergodic joint state process $\mathbf{X}(t)$ and joint action process $\mathbf{A}(t)$, the 
resulting average-reward MARBP is
\begin{equation}
\label{eq:acmarbp}
\maxobj_{\boldsymbol{\pi} \in \boldsymbol{\Pi}^{\textup{SR}}(M)} \Ex^{\boldsymbol{\pi}} \left[\sum_{k \in \mathbb{K}} R_k\big(\widetilde{X}_k, \widetilde{A}_k\big)\right],
\end{equation}
where $\widetilde{\mathbf{X}} = (\widetilde{X}_k)$ and $ \widetilde{\mathbf{A}} = (\widetilde{A}_k)$ are random vectors with 
the corresponding steady-state distributions.

Given the intractability of such a problem, 
\citet{whit88} defined  an index for an RB project, and proposed as a heuristic the resulting index policy:  at each time, engage a subset of projects with the $M$ currently highest  index values.

\begin{figure}[!ht]
\centering
\includegraphics[height=2in]{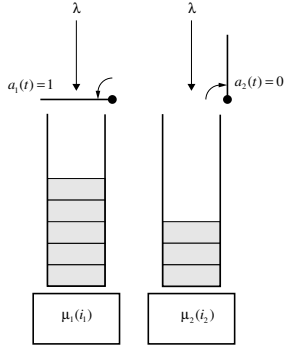}
\caption{Routing to two queues as a two-armed restless bandit.}
\label{fig:marbpf}
\end{figure}

To formulate (\ref{eq:etmp2}) as an MARBP,   we use the idea in \citet[Sec.\ 8.1]{nmmp02} of viewing each queue $k$ as being fed with  its own independent Poisson arrival stream with rate $\lambda$,  subject to \emph{dynamic control of admission}, 
identifying each controlled queue with an RB project.  
One can visualize the admission control actions at each queue as being carried out by a \emph{dedicated gatekeeper}, who 
can open or shut the \emph{entry gate} to his queue. 
We thus interpret the active action $A_k(t) = 1$ for the $k$th RB   as shutting the entry gate to queue $k$, and the passive action $A_k(t) = 0$ as opening it. See
Figure \ref{fig:marbpf}.

The optimal routing problem (\ref{eq:etmp2}) is formulated as a variant of the MARBP in
 (\ref{eq:acmarbp}),  
differing
 from the latter in two respects: (i) whereas in  (\ref{eq:acmarbp}) both actions are available at  
 every state, we  incorporate the requirement  (\ref{eq:ak1}) that the active action be taken 
in project $k$ at state $x_k = n_k$ (the entry gate to a full queue must be closed); and (ii) we replace the requirement that $M \triangleq K-1$ projects be active at each time by the weaker
constraint  (\ref{eq:spconstr}) that $K-1$ projects be active when the system is not full; when the system is full,  (i) implies that the $K$ projects will be active.

\subsection{Relaxed Problem, Lagrangian Relaxation and Decomposition, and a Bound}
\label{s:rplrd}
Recall that  $\boldsymbol{\Pi}^{\textup{SR}}$ denotes the class of admissible policies for problem (\ref{eq:etmp2}). 
Denote by $\widehat{\boldsymbol{\Pi}}^{\textup{SR}}$  
the wider class of policies that is obtained from $\boldsymbol{\Pi}^{\textup{SR}}$ by dropping the requirement in (\ref{eq:spconstr}) that exactly $K-1$ gates be shut when some queue is not full.
We now replace the sample-path constraint (\ref{eq:spconstr})  with the weaker constraint that the average number of closed gates be at least $K-1$, 
\begin{equation}
\label{eq:avconstr}
\Ex^{\boldsymbol{\pi}}\Bigg[\sum_{k \in \mathbb{K}} \widetilde{A}_k\Bigg] \geqslant K-1,
\end{equation}
to obtain the following constrained-MDP \emph{relaxation} of problem (\ref{eq:etmp2}):
\begin{equation}
\label{eq:rp}
\maxobj_{\boldsymbol{\pi} \in \widehat{\boldsymbol{\Pi}}^{\textup{SR}}, (\ref{eq:avconstr})} \! \Ex^{\boldsymbol{\pi}}\Bigg[\sum_{k \in \mathbb{K}} \bar{\mu}_k\big(\widetilde{X}_k\big)\Bigg].
\end{equation}
We will denote by $V^{\textup{R}}$ the maximum value of problem $(\ref{eq:rp})$. Recall that $V^*$ is the maximum throughput rate in problem $(\ref{eq:etmp2})$.

\begin{lemma}
\label{lma:vrubmtvstar}
$V^{\textup{R}}$  is an upper bound on the maximum throughput $V^*$, i.e., 
$V^* \leqslant V^{\textup{R}}$.
\end{lemma}
\begin{proof}
Any feasible policy $\boldsymbol{\pi} \in \widehat{\boldsymbol{\Pi}}^{\textup{SR}}$ for $(\ref{eq:etmp2})$ is also feasible for relaxed problem $(\ref{eq:rp})$, since it belongs to $\widehat{\boldsymbol{\Pi}}^{\textup{SR}}$ and satisfies (\ref{eq:avconstr}). Since the optimization objective is the same in both problems,  their optimal values must satisfy that $V^* \leqslant V^{\textup{R}}$.
\end{proof}

Unlike the original problem (\ref{eq:etmp2}), the \emph{relaxed problem}  (\ref{eq:rp}) will turn out to be tractable. To address it we dualize the constraint (\ref{eq:avconstr}), scaled by $\lambda$,
with a multiplier $\eta \geqslant 0$, which yields  the   \emph{Lagrangian relaxation} 
\begin{equation}
\label{eq:lr}
\begin{split}
\maxobj_{\boldsymbol{\pi} \in \widehat{\boldsymbol{\Pi}}^{\textup{SR}}} \,  \Ex^{\boldsymbol{\pi}}\Bigg[\sum_{k \in \mathbb{K}} \bar{\mu}_k\big(\widetilde{X}_k\big) 
  - \eta \lambda
  \Big\{K-1  - \sum_{k \in \mathbb{K}}   \widetilde{A}_k\Big\}\bigg],
\end{split}
\end{equation}
whose
optimal value satisfies that $V^{\textup{L}}(\eta) \geqslant V^{\textup{R}}$.
Independence of the queues' state transitions allows us to
\emph{decompose} (\ref{eq:lr}) into the $K$ individual queue
 subproblems
\begin{equation}
\label{eq:spsk}
\maxobj_{\pi_k \in \Pi_k^{\textup{SR}}} \! \Ex^{\pi_k}\Big[\bar{\mu}_k\big(\widetilde{X}_k\big) + \eta \lambda \widetilde{A}_k\Big],
\end{equation}
where $\Pi_k^{\textup{SR}}$ is the class of
stationary randomized  admission control policies for operating queue $k$ \emph{in isolation}, which are required to satisfy (\ref{eq:ak1}) and to induce an ergodic state process $\widetilde{X}_k$ and action process $\widetilde{A}_k$. 

Note that (\ref{eq:spsk}) is an \emph{optimal admission control subproblem} for queue $k$, whose objective  balances average throughput rewards and  \emph{rejection subsidies}, 
the latter being earned at rate $\eta$ per rejected job.
Denoting its optimal value  by 
$V_k^*(\eta)$,  the optimal value of Lagrangian relaxation (\ref{eq:lr}) decomposes as
\begin{equation}
\label{eq:optvalrel}
V^{\textup{L}}(\eta) = -(K-1) \lambda \eta
+ \sum_{k \in \mathbb{K}}  V_k^*(\eta).
\end{equation}

Such an approach  
decentralizes the solution to relaxed problem
(\ref{eq:rp}), with $\eta$ playing the role of an economic incentive used by
 the central controller to subsidize
gatekeepers for  rejecting arrivals, while
leaving them to solve
their own subproblems (\ref{eq:spsk}).
A full decentralization  will result by pitching 
$\eta$ to a critical level $\eta^*$  under which
$V^{\textup{L}}(\eta^*) = V^{\textup{R}}$.
Existence of such an $\eta^*$ follows from  \emph{LP strong duality}, drawing on the well-known LP formulation for solving the DP equations of unichain average-reward finite state and action MDPs.
Thus, $\eta^*$ is obtained by solving the  \emph{Lagrangian dual problem}
\begin{equation}
\label{eq:ldp}
\minobj_{\eta \geqslant 0} V^{\textup{L}}(\eta).
\end{equation}

The latter, being a \emph{scalar convex minimization problem}, can be solved efficiently if the 
$V_k^*(\eta)$  can be quickly evaluated. 
Such is the case in the present setting, as the following new result shows. 
Let
\begin{equation}
\label{eq:psumbk}
B(\lambda) \triangleq \sum_{k \in \mathbb{K}} B_{m_k, n_k}(\lambda/\mu_k), 
\end{equation}
where $B_{m_k, n_k}(r_k)$ denotes the blocking
probability in the M/M/$m_k$/$n_k$ queue with offered load $r_k$.

\begin{lemma}
\label{lma:oci}
 \begin{itemize}
 \item[\textup{(a)}] 
It is optimal for $(\ref{eq:spsk})$ to\textup{:} reject all arrivals if $\eta > 1;$ accept arrivals when the queue is not full if $\eta < 1;$
take any action  if $\eta = 1;$
 \item[\textup{(b)}] the optimal value of queue $k$ subproblem $\textup{(\ref{eq:spsk})}$ is given by
 \[
 V_k^*(\eta) = 
 \begin{cases}
 \lambda \big[1 - (1-\eta) B_{m_k, n_k}(\lambda/\mu_k)\big] & \text{ if } \eta < 1; \\
 \lambda  \eta & \text{ if } \eta \geq 1;
 \end{cases}
 \]
 \item[\textup{(c)}] the optimal value of the Lagrangian relaxation \textup{(\ref{eq:lr})}  is given by
 \[
 V^{\textup{L}}(\eta) =
 \begin{cases}
\displaystyle  \lambda \big[K -  B(\lambda) + \big(B(\lambda) - (K-1)\big) \eta \big] & \text{ if } \eta < 1; \\
 \lambda  \eta & \text{ if } \eta \geq 1.
 \end{cases}
 \]
 \end{itemize}
 \end{lemma}
\begin{proof}
 Part (a) follows since admission control subproblem $\textup{(\ref{eq:spsk})}$ balances average throughput rewards, which are earned at a rate of $1$ per completed job, against average rejection subsidies, which are earned at a rate of $\eta$ per rejected job.
Therefore, e.g., in the case $\eta > 1$,  rejecting each arrival and hence collecting the resulting subsidy of $\eta$ is more profitable than accepting it only to receive a lesser reward of $1$ when the job completes. The cases $\eta < 1$ and $\eta = 1$ are argued similarly.
 
 Part (b) follows by evaluating the optimal policies identified in part (a).
Thus, in the case $\eta < 1$,  the optimal policy accepts all arrivals to queue $k$ when it  is not full, inducing an M/M/$m_k$/$n_k$ queue. 
Hence, under such an optimal policy,  arrivals to the queue  are  lost with probability
\[
\Ex\big[\widetilde{A}_k\big] = \Prob\big\{\widetilde{A}_k = 1\big\} = \Prob\big\{\widetilde{X}_k = n_k\big\} = B_{m_k, n_k}(\lambda/\mu_k),
\]
and the queue throughput is 
$\Ex\big[\bar{\mu}_k\big(\widetilde{X}_k\big)\big] = \lambda - \lambda B_{m_k, n_k}(\lambda/\mu_k).$
Therefore, the optimal objective value for subproblem $\textup{(\ref{eq:spsk})}$ is given by
\[
\lambda - \lambda B_{m_k, n_k}(\lambda/\mu_k) + \eta \lambda B_{m_k, n_k}(\lambda/\mu_k) = 
\lambda - (1-\eta) \lambda B_{m_k, n_k}(\lambda/\mu_k).
\]

 In the case $\eta \geq 1$, it is optimal to reject all arrivals to queue $k$, and hence $\widetilde{A}_k \equiv 1$ and $\widetilde{X}_k \equiv 0$. Therefore, $\Ex\big[\widetilde{A}_k\big] = 1$ and $\Ex\big[\bar{\mu}_k\big(\widetilde{X}_k\big)\big] = 0$, which yields that the optimal objective value for subproblem $\textup{(\ref{eq:spsk})}$ takes on the value $ \lambda \eta$.

 Part (c)  follows directly from part (b), (\ref{eq:optvalrel}) and (\ref{eq:psumbk}).
 \end{proof}

We are now ready to solve the Lagrangian dual problem (\ref{eq:ldp}), giving  in closed form its optimal value, which equals that of  relaxed problem (\ref{eq:rp}), and provides an upper bound on the optimal throughput $V^*$ for  problem $(\ref{eq:etmp2})$.
\begin{proposition}[Throughput upper bound]
 \label{pro:ubndthrpt}
 If $B(\lambda) > K-1$ $($resp.\ if $B(\lambda) \leq K-1),$ then $\eta^* = 0$ $($resp.\ $\eta^* = 1)$ solves the Lagrangian dual problem \textup{(\ref{eq:ldp})}, and hence its  optimal value is $V^{\textup{R}} = V^{\textup{L}}(\eta^*) =  \big[K - B(\lambda)\big] \lambda$ $($resp.\ $V^{\textup{R}} = V^{\textup{L}}(\eta^*) = \lambda)$.
\end{proposition}
 \begin{proof}
 The result follows immediately from (\ref{eq:ldp}) and  Lemma \ref{lma:oci}(c) through elementary arguments. 
\end{proof}

Note that the upper bound $V^{\textup{R}}$ on the optimal throughput $V^*$ given in Proposition \ref{pro:ubndthrpt} is nontrivial  only in the case $B(\lambda) > K-1$, in which $V^{\textup{R}} = \big[K - B(\lambda)\big] \lambda < \lambda$.
One sees from (\ref{eq:psumbk}) that the inequality $B(\lambda) > K-1$ will hold when the  arrival rate $\lambda$ is  high enough, since $B(\lambda) \to K$ as $\lambda \to \infty$.

Note further that Proposition \ref{pro:ubndthrpt} also yields a lower bound on the minimum loss rate $J^*$ for problem (\ref{eq:pa2a}): since $J^* + V^* = \lambda$, we obtain 
\begin{equation}
\label{eq:jstlb}
J^* \geqslant  \big[B(\lambda) - (K-1)\big] \lambda.
\end{equation}
 
\subsection{The RB Index does not Yield a Particular Policy}
\label{s:iaip}
We next adapt to  the admission control subproblem (\ref{eq:spsk}) for queue $k$ the concept of indexability  for an RB introduced in \citet{whit88} and extended in \citet{nmmp02}.
Consider the collection of single-queue admission control subproblems (\ref{eq:spsk}), parameterized by the rejection subsidy $\eta$, which is allowed to range over the real line $\mathbb{R}$. We will thus refer to the $\eta$-subsidy subproblem. 
Let us say that such a problem collection is \emph{indexable} if there exists a \emph{break-even rejection subsidy} $\eta_k^{\textup{RB}}(x_k)$ attached to every queue state $0 \leqslant x_k < n_k$, 
such that, for any given rejection subsidy $\eta \in \mathbb{R}$,   it is optimal for the $\eta$-subsidy subproblem to reject an arrival finding $0 \leqslant x_k < n_k$ jobs present  iff 
$\eta_k^{\textup{RB}}(x_k) \leqslant \eta$.
We will refer to $\eta_k^{\textup{RB}}(x_k)$ as the \emph{RB} or \emph{Whittle index} for queue $k$. 

It follows immediately from the definition of indexability and Lemma \ref{lma:oci}(a) that 
 each subproblem (\ref{eq:spsk}) is indexable, with the RB index being constant across states: $\eta_k^{\textup{RB}}(x_k) \equiv 1$ 
We thus find that the RB index for this model   does not yield a particular policy --- since the resulting index policy prescribes  routing an arrival at time $t$ to a nonfull queue, if any,  with currently highest index $\eta_k^{\textup{RB}}\big(X_k(t)\big)$, breaking ties arbitrarily.

\subsection{A Tie-Breaking Second-Order Index Policy}
\label{s:tbsoip}
To break the ties caused by use of the RB index, we will use a second-order RB index, as introduced in \citet[Sec.\ 2.2]{nmqs06} for a different model.
To define the second-order index, consider the $\alpha$-discounted version of queue $k$'s admission control subproblem (\ref{eq:spsk}), where $\alpha > 0$ denotes the discount rate:
\begin{equation}
\label{eq:spskd}
\maxobj_{\pi_k \in \Pi_k^{\textup{SR}}} \! \Ex_{x_k}^{\pi_k}\bigg[\int_0^\infty e^{-\alpha t} \left\{\bar{\mu}_k\big(X_k(t)\big) + \eta \lambda A_k(t)\right\}  \, dt\bigg].
\end{equation}
In (\ref{eq:spskd}), $\Ex_{x_k}^{\pi_k}[\cdot]$ denotes expectation under policy $\pi_k$, conditioned on $X_k(0) = x_k$.

Similarly as for (\ref{eq:spsk}), let us say that subproblem (\ref{eq:spskd}) is indexable if there exists a break-even rejection subsidy $\eta_k^{\textup{RB}}(x_k; \alpha)$ for each queue state $0 \leqslant x_k < n_k$, 
such that, for any rejection subsidy $\eta \in \mathbb{R}$,   it is optimal to reject an arrival finding $0 \leqslant x_k < n_k$ jobs present, regardless of the initial state,  iff 
$\eta_k^{\textup{RB}}(x_k; \alpha) \leqslant \eta$.
We will refer to $\eta_k^{\textup{RB}}(x_k; \alpha)$ as the queue's \emph{$\alpha$-discount RB index}.

Further, if  the first-order Maclaurin series expansion of $\eta_k^{\textup{RB}}(x_k; \alpha)$ is 
\begin{equation}
\label{eq:fomse}
\eta_k^{\textup{RB}}(x_k; \alpha) = 1 - \theta_k^{\textup{RB}}(x_k) \alpha + o(\alpha), \quad \text{as } \alpha \searrow 0,
\end{equation}
then we define the tie-breaking second-order RB index to be $\theta_k^{\textup{RB}}(x_k)$. 
It is the latter that we will use to define a heuristic routing index policy for problem (\ref{eq:etmp2}): route an arrival at time $t$ to a nonfull queue, if any, with \emph{lowest} index $\theta_k^{\textup{RB}}\big(X_k(t)\big)$.

We next turn to showing that subproblem (\ref{eq:spskd}) is indexable, and to evaluating the second-order index $\theta_k^{\textup{RB}}(x_k)$. 
We will attain both goals by reducing (\ref{eq:spskd}) to an equivalent subproblem that, unlike (\ref{eq:spskd}), satisfies the 
sufficient indexability conditions for admission control problems given in \citet[Sec.\ 7]{nmmp02}.
The following invariance identity will play a key role.

\begin{lemma}
\label{lma:accarg} Under any admission policy $\pi_k \in \Pi_k^{\textup{SR}}$ for queue $k,$ 
\[
\Ex_{x_k}^{\pi_k}\left[\int_0^\infty e^{-\alpha t} \big\{\alpha X_k(t) +
    \bar{\mu}_k\big(X_k(t)\big) - \lambda \big(1-A_k(t)\big)\big\}\, dt\right] \equiv
x_k.
\]
\end{lemma}
\begin{proof}
Suppose that linear holding costs are  incurred over time at unit rate. Then, 
the expected total discounted (ETD) holding cost, 
$\Ex_{x_k}^{\pi_k}\Big[\int_0^\infty e^{-\alpha t} X_k(t) \, dt\Big]$, can be equivalently accounted for via job-oriented lump charges
as  in \citet{bell71}. To do so, charge $1/\alpha$ for each 
job initially present and for each 
job admitted thereafter, as though they were to stay
forever in the system, which gives the  ETD charge  $x_k/\alpha
+ (\lambda/\alpha)\Ex_{x_k}^{\pi_k}\Big[\int_0^\infty e^{-\alpha
  t}\big(1-A_k(t)\big) \, dt\Big]$.
Then, refund to each job upon departure the amount $1/\alpha$,
which gives the ETD refund 
$(1/\alpha) \Ex_{x_k}^{\pi_k}\Big[\int_0^\infty e^{-\alpha t} \bar{\mu}_k\big(X(t)\big) \,
dt\Big]$. Equating the expected costs under both charging schemes yields the result. 
\end{proof}

We next draw on Lemma \ref{lma:accarg} to  reformulate the objective of subproblem (\ref{eq:spskd})
in a form that is more convenient for our purposes.

\begin{lemma}
\label{lma:refspspd} The objective of subproblem \textup{(\ref{eq:spskd})} can be
reformulated as
\[
x_k + \lambda/\alpha - \Ex_{x_k}^{\pi_k}\left[\int_0^\infty e^{-\alpha t} \big\{\alpha X_k(t)  + (1-\eta) \lambda A_k(t)\big\}\, dt\right].
\]
\end{lemma}
\begin{proof}
By Lemma \ref{lma:accarg}, we have
\[
\Ex_{x_k}^{\pi_k}\left[\int_0^\infty e^{-\alpha t} \bar{\mu}_k\big(X_k(t)\big)\, dt\right]  
=
x_k - \Ex_{x_k}^{\pi_k}\left[\int_0^\infty e^{-\alpha t} \big\{\alpha X_k(t)  - \lambda \big(1-A_k(t)\big)\big\}\, dt\right],
\]
and hence we can write the objective of subproblem (\ref{eq:spskd}),
\[
\Ex_{x_k}^{\pi_k}\bigg[\int_0^\infty e^{-\alpha t} \left\{\bar{\mu}_k\big(X_k(t)\big) + \eta \lambda A_k(t)\right\}  \, dt\bigg],
\]
as the given expression.
\end{proof}

Consider now the discounted-cost subproblem
\begin{equation}
\label{eq:minhcrc}
\minobj_{\pi_k \in \Pi_k^{\textup{SR}}} \, \Ex_{x_k}^{\pi_k}\left[\int_0^\infty e^{-\alpha t} \big\{X_k(t)  + \nu \lambda A_k(t)\big\}\, dt\right],
\end{equation}
which is to find an admission control policy for queue $k$ minimizing the expected discounted value of holding costs and \emph{rejection charges}, 
with the former continuously accruing at unit rate per unit time a job spends in the system, and the latter being incurred at rate $\nu$ per rejected job.
This is a special case of the general optimal admission control problem solved in \citet[Sec.\ 7]{nmmp02} via RB indexation.
The results there ensure existence of the RB  index $\nu_k^{\textup{RB}}(x_k; \alpha)$ attached to  states $0 \leqslant x_k < n_k$, which is such that, for any 
value of the rejection charge $\nu \in \mathbb{R}$, it is optimal for problem (\ref{eq:minhcrc}) to reject an arrival in state $x_k$ iff $\nu_k^{\textup{RB}}(x_k; \alpha) \geqslant \nu$.

We draw on \citet[Sec.\ 7]{nmmp02} to address the average-cost subproblem
\begin{equation}
\label{eq:minhcrac}
\minobj_{\pi_k \in \Pi_k^{\textup{SR}}} \, \Ex^{\pi_k}\left[\int_0^\infty \big\{\widetilde{X}_k  + \nu \lambda \widetilde{A}_k\big\}\, dt\right],
\end{equation}
as results there ensure  existence of the RB index $\nu_k^{\textup{RB}}(x_k) = \lim_{\alpha \searrow 0} \nu_k^{\textup{RB}}(x_k; \alpha)$, which is  such that, for any 
$\nu \in \mathbb{R}$, it is optimal for (\ref{eq:minhcrac}) to reject an arrival in state $0 \leqslant x_k< n_k$ iff $\nu_k^{\textup{RB}}(x_k) \geqslant \nu$.

We are now ready to establish that subproblem (\ref{eq:spskd}) is indexable, and to evaluate the second-order index $\theta_k^{\textup{RB}}(x_k)$ defined via (\ref{eq:fomse}). 

\begin{proposition}
\label{pro:indx}
\begin{itemize}
\item[\textup{(a)}] Subproblem $(\ref{eq:spskd})$
  is indexable, with its RB index $\eta_k^{\textup{RB}}(x_k; \alpha)$ being related to that of subproblem $(\ref{eq:minhcrc})$ by 
  $\eta_k^{\textup{RB}}(x_k; \alpha) = 1 - \alpha \nu_k^{\textup{RB}}(x_k; \alpha);$
\item[\textup{(b)}] we have the Maclaurin expansion $\eta_k^{\textup{RB}}(x_k; \alpha) = 1 - \nu_k^{\textup{RB}}(x_k) \alpha + o(\alpha)$ as $\alpha \searrow 0$, and hence the second-order RB index defined via \textup{(\ref{eq:fomse})} is  $\theta_k^{\textup{RB}}(x_k) = \nu_k^{\textup{RB}}(x_k)$.
\end{itemize}
\end{proposition}
\begin{proof}
(a) It follows from Lemma \ref{lma:refspspd} that subproblem $(\ref{eq:spskd})$ is equivalent to 
\[
\minobj_{\pi_k \in \Pi_k^{\textup{SR}}} \, \Ex_{x_k}^{\pi_k}\left[\int_0^\infty e^{-\alpha t} \big\{X_k(t)  + \frac{1-\eta}{\alpha} \lambda A_k(t)\big\}\, dt\right],
\]
and hence the results in \citet[Sec.\ 7]{nmmp02} ensure that it is optimal to reject an arrival in state $0 \leqslant x_k < n_k$ iff $\nu_k^{\textup{RB}}(x_k; \alpha) \geqslant (1-\eta)/\alpha$, i.e., iff 
$1-\alpha \nu_k^{\textup{RB}}(x_k; \alpha) \leqslant \eta$. This shows that subproblem  $(\ref{eq:spskd})$ is indexable, with RB index $\eta_k^{\textup{RB}}(x_k; \alpha) = 1 - \alpha 
\nu_k^{\textup{RB}}(x_k; \alpha)$.

(b)
This part follows from (a) and the result referred to above on the existence of the RB index $\nu_k^{\textup{RB}}(x_k) = \lim_{\alpha \searrow 0} \nu_k^{\textup{RB}}(x_k; \alpha)$ for average-cost subproblem (\ref{eq:minhcrac}).
\end{proof}

\subsection{Efficient Computation and a Closed Formula for the Index $\theta_k^{\textup{RB}}(x_k)$}
\label{s:eerirb}
In light of the identity $\theta_k^{\textup{RB}}(x_k) \equiv \nu_k^{\textup{RB}}(x_k)$ in Proposition \ref{pro:indx}(b), to evaluate $\theta_k^{\textup{RB}}(x_k)$ we can use  the \emph{coupled first-order recursions} given in \citet[pp.\ 396--397]{nmmp02} for computing the index  $\nu_k^{\textup{RB}}(x_k)$. Besides $\theta_k^{\textup{RB}}(x_k)$, such recursions evaluate quantities $y_k(x_k)$ and $w_k(x_k)$, as follows, where we write 
$\bar{\mu}_k(x_k) \triangleq \min(x_k, m_k) \mu_k$, $\rho_k(x_k) \triangleq
\lambda/\bar{\mu}_k({x_k+1})$, and  $\Delta \bar{\mu}_k(x_k) \triangleq \bar{\mu}_k(x_k) - \bar{\mu}_k(x_k-1)$:  
for $1 \leq x_k < n_k$, 
\begin{equation}
\label{eq:windexqs}
\begin{split}
\theta_k^{\textup{RB}}(x_k)  & = 
\theta_k^{\textup{RB}}(x_k-1) + \frac{\displaystyle{1 - 
  \theta_k^{\textup{RB}}(x_k-1)} \Delta \bar{\mu}_k(x_k+1)}{\displaystyle{\Delta \bar{\mu}_k(x_k+1) + 
  w_k(x_k-1)/\rho_k(x_k-1)}}\,  \\
y_k(x_k) & =  1-\frac{\lambda  \bar{\mu}_k({x_k})/y_k({x_k-1})}{\big[\lambda +
\bar{\mu}_k({x_k})\big]  \big[\lambda + \bar{\mu}_k(x_k+1)\big]} \\
w_k({x_k})  & = \lambda
\frac{\displaystyle \Delta \bar{\mu}_k(x_k+1) + 
  w_k({x_k-1})/\rho_k({x_k-1})}{y_k(x_k) \big[\lambda + \bar{\mu}_k(x_k+1)\big]},
\end{split}
\end{equation}
with initial values $\theta_k^{\textup{RB}}(0) = 1/\mu_k$, $y_k(0) = 1$, and $w_k(0) = 
\lambda \mu_k/(\lambda + \mu_k)$.

Using (\ref{eq:windexqs}) the  index values $\theta_k^{\textup{RB}}(0), \ldots,
\theta_k^{\textup{RB}}(n_k-1)$ are efficiently evaluated  in $O(n_k)$ time.

We next develop a closed formula for the index $\theta_k^{\textup{RB}}(x_k)$,  by exploiting the representation in \citet[The.\ 6.4(b)]{nmmp02} for $\nu_k^{\textup{RB}}(x_k)$, given by
\begin{equation}
\label{eq:cfmpim}
\theta_k^{\textup{RB}}(x_k) = \frac{1}{\lambda}
\frac{L_{m_k, x_k+1}(r_k) - L_{m_k, x_k}(r_k)}{B_{m_k, x_k}(r_k) -
  B_{m_k, x_k+1}(r_k)},
\end{equation} 
where $L_{m_k, x_k}(r_k)$ and $B_{m_k, x_k}(r_k)$ denote the mean number in system and the blocking
probability for the M/M/$m_k$/$x_k$ queue with offered load $r_k = \lambda/\mu_k$, respectively.
The following result refers to the Erlang-B formula $B_{m_k}(r_k)$ for the blocking probability in the M/M/$m_k$/$m_k$ queue with offered load $r_k$, and to the Erlang-C formula $C_{m_k}(r_k)$, which is used even in the case $r_k \geq m_k$, where it does not have the interpretation for the case  $r_k < m_k$ as the  
 delay probability for the M/M/$m_k$ queue with offered load $r_k$.
We further denote by $C_{m_k}'(r_k)$ the derivative of $C_{m_k}(r_k)$ with respect to $r_k$, and by $\rho_k \triangleq r_k/m_k$ the offered load per server.

\begin{proposition}
\label{pro:mpifm} The RB index $\theta_k^{\textup{RB}}(x_k)$ has the following evaluation$:$
$\theta_k^{\textup{RB}}(x_k) = 1/\mu_k$ for $0 \leq x_k < m_k;$ and, for $m_k \leq x_k < n_k,$
\[
\theta_k^{\textup{RB}}(x_k) = 
\displaystyle \bigg[ \frac{  \rho_k C_{m_k}(r_k) (\rho_k^{x_k-m_k+1} - 1)}{m_k \mu_k (\rho_k-1)^2}-\frac{x_k+1-r_k}{m_k \mu_k (\rho_k-1)}\bigg], \, \text{if } \rho_k \neq 1,
\]
and
\[
\theta_k^{\textup{RB}}(x_k) = 
\displaystyle \frac{\frac{1}{2}\big(x_k-m_k + 2 + 2 m_k  C_{m_k}'(m_k)\big)(x_k-m_k+1) + m_k}{m_k \mu_k},\, \text{if } \rho_k = 1.
\]
\end{proposition}
\begin{proof}
Two cases need be distinguished. 
If  $x_k < m_k$, the M/M/$m_k$/$x_k$ queue 
reduces to the 
M/M/$x_k$/$x_k$ queue. Hence, $B_{m_k, x_k}(r_k) = B_{x_k}(r_k)$, 
$L_{m_k,  x_k}(r_k) =  r_k 
\big[1-B_{x_k}(r_k)\big]$ and, therefore,
$\theta_k^{\textup{RB}}(x_k) = 1/\mu_k$.
In the second case, when $m_k \leq x_k < n_k$, substituting 
for $B_{m_k,  x_k}(r_k)$ and $L_{m_k, x_k}(r_k)$  their closed formulae (see, e.g., \citet[Sec.\ 2.5]{grossetal08}), and
simplifying  the resulting expressions, gives the result when $\rho_k \neq 1$. 

The formula for the case $\rho_k = 1$ is obtained by replacing (in the formula for the case $\rho_k \neq 1$) $C_{m_k}(r_k)$ with its first-order Taylor expansion $1 - m_k C_{m_k}'(m_k) (1-\rho_k)$, and simplifying the resulting expression as $\rho_k \to 1$.
\end{proof}

To evaluate the required $C_{m_k}'(m_k)$ when $\rho_k = 1$, we can draw on the identities
\[
C_{m_k}(r_k) = \frac{m_k B_{m_k}(r_k)}{m_k - r_k + r_k B_{m_k}(r_k)},  \,
B_{m_k}'(r_k) = \frac{m_k - r_k + r_k B_{m_k}(r_k)}{r_k} B_{m_k}(r_k),
\]
where $B_{m_k}'(r_k)$ is the derivative of $B_{m_k}(r_k)$ with respect to $r_k$, which yield that
$B_{m_k}'(m_k) = B_{m_k}(m_k)^2$ and $C_{m_k}'(m_k) = \big[1 - B_{m_k}(m_k)\big]/\big(m_k B_{m_k}(m_k)\big)$.
 
 \section{PI Index Method}
\label{s:obspim}
This section discusses application of the PI index method, based on performing one PI step over an OBS, extending the results in \citet{krish88}.

\subsection{First Stage: OBS Method}
\label{s:fsobsm}
This section considers the OBS method to design a static routing policy for the present model, extending work in 
 \citet{yaoshant87}. Whereas the latter paper addressed only the special no-waiting-room case $m_k \equiv n_k$, the following discussion applies to the general case $m_k \leqslant n_k$.
 
\subsubsection{BS Policies.}
\label{s:bsp}
Consider a \emph{Bernoulli splitting} (BS) where each arrival is independently routed to queue $k$ with a fixed probability $p_k$.
Note that BS policies can violate requirement (\ref{eq:ak1}) by routing an arrival to a full queue, and hence they do not belong to the class $\boldsymbol{\Pi}^{\textup{SR}}$ of admissible
policies considered in Section \ref{s:mod}.

To accommodate such static policies, we  reformulate problem (\ref{eq:pa2a}) as 
\begin{equation}
\label{eq:pa2}
\minobj_{\boldsymbol{\pi} \in \widetilde{\boldsymbol{\Pi}}^{\textup{SR}}} 
\lambda \Ex^{\boldsymbol{\pi}}\left[\sum_{k \in \mathbb{K}}  \big(1-\widetilde{A}_k\big)  1_{\{\widetilde{X}_k = n_k\}}\right],
\end{equation}
where 
$\widetilde{\boldsymbol{\Pi}}^{\textup{SR}}$ is the class of stationary randomized policies obtained from $\boldsymbol{\Pi}^{\textup{SR}}$ by allowing randomized actions, 
dropping  the constraint (\ref{eq:ak1}), and replacing (\ref{eq:spconstr}) with 
 \begin{equation}
 \label{eq:spconstr2}
 \sum_{k \in \mathbb{K}} A_k(t)  = K-1,  \quad t \geqslant 0.
 \end{equation}
 Thus, under policies in $\widetilde{\boldsymbol{\Pi}}^{\textup{SR}}$ every arrival is routed to some queue, which can be full. 
 Note that the objective in (\ref{eq:pa2}) is the (long-run) expected proportion of jobs that are lost due to being routed to a full queue, and that
 problems (\ref{eq:pa2a}) and (\ref{eq:pa2}) are equivalent, provided the following result holds.
 
 \begin{conjecture}
 \label{con:probequiv}
An optimal policy for problem \textup{(\ref{eq:pa2})} routes an arrival to a nonfull queue, if any is available. 
 \end{conjecture}  
While Conjecture \ref{con:probequiv} appears intuitively clear, the author has not found a proof in the literature. 

Under a BS, the queues are decoupled, with queue $k$ behaving as an M/M/$m_k$/$n_k$ queue with arrival rate 
$\lambda_k \triangleq \lambda p_k$.
Denoting by  $r_k(\lambda_k) \triangleq 
\lambda_k/ \mu_k$ the \emph{offered load to queue} $k$,
the total  loss rate is
$\sum_{k \in \mathbb{K}} \phi_k(\lambda_k)$,
where $\phi_k(\lambda_k) \triangleq \lambda_k B_{m_k, n_k}\big(r_k(\lambda_k)\big)$
and $B_{m, n}(r)$ denotes
the  blocking probability in the
M/M/$m$/$n$ queue with offered load $r$.
Note that the derivative of the loss rate $\phi_k(\lambda_k)$ for queue $k$ has the evaluation
\begin{equation}
\label{eq:fkprime}
\phi_k'(\lambda_k) = B_{m_k, n_k}\big(r_k(\lambda_k)\big) +
r_k(\lambda_k) B_{m_k, n_k}'\big(r_k(\lambda_k)\big),
\end{equation}
from which it follows that $\phi_k'(0) = 0$, a result we will use below.

Each $\phi_k(\lambda_k)$ is 
continuously differentiable. \citet[Cor.\ 5.i]{shanthiYao92} shows that $\phi_k(\lambda_k)$ is 
   convex nondecreasing.
\citet[Th.\ 2.9 (2.26i)]{liyanShanthik92} states the stronger result
 (though referring to the former paper for a proof) that  $\phi_k(\lambda_k)$ is 
  increasing and strictly convex. 

\begin{proposition}[Liyanage, Shanthikumar, and Yao]
\label{con:fklambcnvxinc}
$\phi_k(\lambda_k)$ is increasing
and
 strictly convex.
\end{proposition}

\subsubsection{Computing the OBS}
\label{s:cobs}
The OBS method is to find a BS $\boldsymbol{\lambda}^* = (\lambda_k^*)$ that 
solves the  NLP
\begin{equation}
\label{eq:nlp}
\minobj \, \Bigg\{\sum_{k \in \mathbb{K}} \phi_k(\lambda_k)\colon 
\sum_{k \in \mathbb{K}} \lambda_k = \lambda,
 0 \leqslant \lambda_k \leqslant \lambda, \, k \in \mathbb{K}\Bigg\}.
\end{equation}
Although the constraints $\lambda_k \leqslant \lambda$ in (\ref{eq:nlp}) are redundant,  we include them for algorithmic reasons, as discussed below.

Problem  (\ref{eq:nlp}) is an example of the much studied separable, convex, and differentiable nonlinear resource allocation problem. Nevertheless,  
 we analyze it next to exploit its special properties.

Proposition \ref{con:fklambcnvxinc}, together with the compactness of  the  feasible solutions region, ensures existence of a unique $\boldsymbol{\lambda}^*$ attaining
the
minimum cost $J^{\textup{BS}}$ of (\ref{eq:nlp}), which is characterized as the unique feasible solution satisfying the  \emph{first-order Karush--Kuhn--Tucker} (KKT) \emph{conditions}: There exists a
  unique multiplier  $y^* \in \mathbb{R}$ attached to the equality
constraint  in (\ref{eq:nlp}) such that
\begin{equation}
\label{eq:kktc}
\begin{split}
\phi_k'(\lambda_k^*)  & \geqslant y^* \text{ if } \lambda_k^* =
0;
\phi_k'(\lambda_k^*)   = y^* \text{ if } 0 < \lambda_k^* <
\lambda; 
\phi_k'(\lambda_k^*)  \leqslant y^* \text{ if } \lambda_k^* =
\lambda.
\end{split}
\end{equation}

To find the optimal $(\boldsymbol{\lambda}^*, y^*)$ pair, we consider the \emph{Lagrangian function}
\begin{equation}
\label{eq:ly}
L^{\textup{BS}}(y) \triangleq \min_{\substack{\lambda_k \in [0, \lambda] \\ k \in \mathbb{K}}}
 \sum_{k \in \mathbb{K}} \phi_k(\lambda_k) + y \Bigg[\lambda - \sum_{k \in \mathbb{K}} \lambda_k\Bigg] = 
 \lambda y  + \sum_{k \in \mathbb{K}} \min_{\lambda_k \in [0, \lambda]} \big[\phi_k(\lambda_k) - \lambda_k y\big],
\end{equation}
which is concave in $y$ and gives a lower bound $L^{\textup{BS}}(y) \leqslant J^{\textup{BS}}$. 
Since the right-hand side in (\ref{eq:ly}) decomposes  into $K$ one-variable strictly convex subproblems, 
we can evaluate $L^{\textup{BS}}(y)$ by finding the unique optimal solution $\lambda_k^*(y)$ to each one, so 
\begin{equation}
\label{eq:lyev}
L^{\textup{BS}}(y)  = \lambda y  + \sum_{k \in \mathbb{K}} \big[\phi_k\big(\lambda_k^*(y)\big) - \lambda_k^*(y) y\big].
\end{equation}

The following result is immediate.
\begin{lemma}
\label{lma:lambdakstar}
\begin{itemize}
\item[\textup{(a)}]
$\lambda_k^*(y) = 0$ 
if $y \leqslant \phi_k'(0)$ $($i.e., if $y \leqslant 0$, since $\phi_k'(0) = 0);$
\item[\textup{(b)}] $\lambda_k^*(y) = \lambda$  if $y \geqslant \phi_k'(\lambda);$ 
\item[\textup{(c)}] if
$0 < y  < \phi_k'(\lambda)$, $\lambda_k^*(y)$ is the unique root in $(0, \lambda)$ of $\phi_k'(\lambda_k) = y;$ 
\item[\textup{(d)}] 
$\lambda_k^*(y)$ is continuous nondecreasing in $y \in \mathbb{R}$, being increasing over $y \in [0, \phi_k'(\lambda)].$
Hence, $\sum_{k \in \mathbb{K}} \lambda_k^*(y)$ is continuous increasing over 
$y \in [0, \max_{k \in \mathbb{K}} \phi_k'(\lambda)].$
\end{itemize}
\end{lemma}

The optimal $y^*$ is found by solving the \emph{Lagrangian dual problem}
\begin{equation}
\label{eq:ldp1}
\maxobj \big\{L^{\textup{BS}}(y)\colon y \in \mathbb{R}\big\},
\end{equation}
which, due to Proposition \ref{con:fklambcnvxinc}, has the same optimal value $J^{\textup{BS}}$ as (\ref{eq:nlp}).

The next result, which exploits the properties discussed above of the functions $\phi_k(\lambda_k)$, gives the key to finding the solutions to dual  and primal problems  (\ref{eq:ldp1}) and (\ref{eq:nlp}). Let $b \triangleq \min_{k \in \mathbb{K}} \phi_k'(\lambda)$.
Part (a) identifies a finite open interval containing the optimal $y^*$, while  part (b) ensures that the OBS routes traffic to every queue.

\begin{proposition}
\label{pro:ldp}
\begin{itemize}
\item[\textup{(a)}] Lagrangian dual problem \textup{(\ref{eq:ldp1})} is equivalent to 
\begin{equation}
\label{eq:ldp2}
\maxobj\big\{L^{\textup{BS}}(y)\colon y \in (0, b)\big\},
\end{equation}
with its optimal solution $y^*$ being the unique root in $(0, b)$ of the equation 
\begin{equation}
\label{eq:ystareq}
\sum_{k \in \mathbb{K}} \lambda_k^*(y) = \lambda;
\end{equation}
\item[\textup{(b)}]  the OBS of the arrival stream is $\boldsymbol{\lambda}^* \triangleq \big(\lambda_k^*(y^*)\big)$, which satisfies that
\begin{equation}
\label{eq:lambdapos}
0 < \lambda_k^* < \lambda, \quad k \in \mathbb{K}.
\end{equation}
\end{itemize}
\end{proposition}
\begin{proof}
The proofs of both parts are  intertwined.
We start by showing that two of the three cases in the KKT conditions (\ref{eq:kktc}) can be dismissed. 
Thus, if the OBS had $\lambda_k^* = 0$ for some station $k$, then the multiplier would be $y^* \leq \phi_k'(0) = 0$. But then, for any station $l$ with $\lambda_l^* > 0$ we would have that $y^* \geq \phi_l'(\lambda_l^*) > 0$, a contradiction, which discards the case $\lambda_k^* = 0$.
And, if the OBS had $\lambda_k^* = \lambda$ for some station $k$, we would then have that $\lambda_l^* = 0$ for some other station $l$, which, as we have just seen, cannot happen. Hence, the only possible case  is $0 < \lambda_k^* < \lambda$, which must hold for every $k$.

Now, it follows from the KKT conditions that it must be $\phi_k'(\lambda_k^*)   = y^*$, and hence $\lambda_k^* = \Lambda_k(y^*)$ for every  $k$. And, finally, since $\sum_k \lambda_k^* = \lambda$, $y^*$ is characterized by the equation (\ref{eq:ystareq}), which has a unique root in $(0, b)$, since the $\Lambda_k(y)$ are increasing in $y$ and satisfy that $\Lambda_k(0) = 0$ and $\Lambda_k(b_k) = \lambda$.
\end{proof}

We can approximate $y^*$ within the desired accuracy using, e.g.,  the  \emph{bisection method} on (\ref{eq:ystareq}), which further yields an
 approximation to $\boldsymbol{\lambda}^*$.

The next result gives the OBS in closed form, assuming equal buffer ($n_k \equiv n$) and server-pool ($m_k \equiv m$)
sizes at the stations, in which case they are shown to have equal offered loads $\lambda_k^*/\mu_k \equiv r^*$.
It extends the result in \citet[Cor.\ 4.1]{yaoshant87} for the case $m_k = n_k \equiv m$ and $\mu_k \equiv \mu$.

\begin{proposition}
\label{pro:osallocm3} 
If $n_k \equiv n$ and  $m_k \equiv m$, then under the OBS all queues have equal offered loads, i.e.,   $\lambda_k^* = \lambda
\mu_k / \sum_{l \in \mathbb{K}}^K \mu_l$, for $k \in \mathbb{K}$.
\end{proposition}
\begin{proof}
From (\ref{eq:fkprime}), 
we have that 
$\phi_k'(\lambda_k) = \phi_{m, n}\big(r_k(\lambda_k)\big)$, with $\phi_{m, n}(r)$ being increasing in $r$ (see Proposition \ref{con:fklambcnvxinc}). 
Hence, to find the root  $\lambda_k^*$ of $\phi_k'(\lambda_k) = y^*$, we can proceed as follows.
Let $r^*$ be the unique root of $\phi_{m, n}(r) = y^*$.
Then, from $r_k(\lambda_k^*) = \lambda_k^*/\mu_k = r^*$ we obtain 
$\lambda_k^* = \mu_k r^*$ for each $k$.
Further, from $\sum_l \lambda_l^* = \lambda$ we obtain 
$r^* \sum_l \mu_l = \lambda$, and hence $r^* = \lambda / \sum_l
\mu_l$ and $\lambda_k^* = \lambda \mu_k / \sum_l
\mu_l$ for each $k$. 
\end{proof}

\subsubsection{Efficient Numerical Evaluations.}
\label{s:ene}
Since  implementation of the BS method requires multiple evaluations of $\phi_k(\cdot)$ and $\phi_k'(\cdot)$ for each queue $k$, it is essential that
such functions be computed as efficiently as possible.
In light of the definition of $\phi_k(\lambda_k)$, and of  (\ref{eq:fkprime}),
the key is to evaluate
 the blocking probability $B_{m, n}(r)$ in the M/M/$m$/$n$
queue with offered load $r$, as well as its derivative with respect to $r$, $B_{m, n}'(r)$.

We first formulate $B_{m, n}(r)$ in terms of the Erlang-B function $B_m(r)$.
Standard results give that, denoting by 
 $\rho \triangleq r/m$ the \emph{offered load per server}  for the queue, 
\begin{equation}
\label{eq:fklambmth} 
B_{m, n}(r) = 
\begin{cases}
\displaystyle \frac{\rho^{n - m} (1-\rho)B_m(r)}{1-\rho + \rho (1-\rho^{n-m}) B_m(r)}, & r \neq m \\
\displaystyle \frac{B_m(m)}{1 + (n - m) B_m(m)}, & r = m.
\end{cases}
\end{equation}

As for $B_m(r)$ and its derivative $B_m'(r)$, by standard results we have
\begin{equation}
\label{eq:bmrrecd}
B_0(r) = 1,  \, B_m(r) = \frac{r B_{m-1}(r)}{m + r B_{m-1}(r)},  \, B_m'(r) = \frac{m - r + r B_{m}(r)}{r} B_{m}(r), \,  m \geqslant 1.
\end{equation}

We next give recursions for computing $B_{m, m}(r)$ and $B'_{m, n}(r)$. 
From (\ref{eq:fklambmth}), 
\begin{equation}
\label{pro:bmnrd1}
B_{m, m}(r) = B_m(r), \, B_{m, n}(r) =  \frac{r B_{m, n-1}(r)}{m + r B_{m, n-1}(r)}, \quad n \geqslant m+1.
\end{equation}
Further, from (\ref{eq:fklambmth}) and (\ref{eq:bmrrecd}) we obtain
\begin{equation}
\label{pro:bmnrd2}
B'_{m, n}(r) = 
\left[\frac{n-r}{\rho} + \bigg(n-r - 
 \frac{1-\rho^{n-m}}{(1-\rho) \rho^{n-m}}\bigg) \frac{B_{m, n}(r)}{1-\rho}
\right] \frac{B_{m, n}(r)}{m}, \, r \neq m,
\end{equation}
\begin{equation}
\label{pro:bmnrd3}
B'_{m, n}(m) = \left[n-m + \big(1+m+n - (1+n-m)^2\big) \frac{B_{m, n}(m)}{2}
\right] \frac{B_{m, n}(m)}{m}.
\end{equation}
The above recursions allow us to compute $B_{m, n}(r)$ and $B'_{m, n}(r)$ in $O(n)$ time. 

\subsubsection{Second Stage: PI Step and Index Policy.}
\label{s:mod2pi}
\citet{krish88} first applied the PI method to the present model, yet focusing only on its second stage (deriving the 
PI index rule), once the first stage (obtaining the OBS) has been carried out. We next briefly reformulate Krishnan's results in the present notation, and further extend them to obtain a new result.  

\citet{krish88} argues from first principles that the PI routing index  has the following evaluation, where $r_k^* \triangleq
\lambda_k^*/\mu_k$, $\rho_k^* \triangleq r_k^*/m_k$, and $\boldsymbol{\lambda}^* = (\lambda_k^*)$ is the OBS:
\begin{equation}
\label{eq:mipii1}
\theta_k^{\textup{PI}}(x_k) = 
\begin{cases}
B_{m_k,  n_k}(r_k^*)/B_{x_k}(r_k^*),  &  0 \leqslant x_k \leqslant m_k, \\
\displaystyle B_{m_k,  n_k}(r_k^*) \frac{1 - \rho_k^* + \rho_k^* (1 -
  (\rho_k^*)^{x_k-m_k}) B_{m_k}(r_k^*)}{(\rho_k^*)^{x_k-m_k} (1-\rho_k^*) B_{m_k}(r_k^*)
},  &  m_k < x_k < n_k.
\end{cases}
\end{equation}
Note that, in the case $\rho_k^* = 1$, the latter formula must be replaced with
\begin{equation}
\label{eq:mipii2b}
\theta_k^{\textup{PI}}(x_k) =  B_{m_k,  n_k}(m_k) \big[x_k-m_k + 1/B_{m_k}(m_k)\big], \quad m_k < x_k < n_k.
\end{equation}
Such expressions arise as solutions to the following  first-order linear recursion:
\begin{equation}
\label{eq:pe1}
\begin{split}
\phi_k^* - \lambda_k^* \theta_k^{\textup{PI}}(0)  & = 0; \quad
\phi_k^* - \lambda_k^* \theta_k^{\textup{PI}}(x_k) + \bar{\mu}_k(x_k)
\theta_k^{\textup{PI}}(x_k-1)  = 0, \quad  1 \leqslant x_k < n_k,
\end{split}
\end{equation}
where $\phi_k^* \triangleq \phi_k\big(\lambda_k^*\big) = \lambda_k^* B_{m_k,  n_k}(r_k^*)$. Thus, (\ref{eq:pe1}) allows to compute the $n_k$ PI index values
$\{\theta_k^{\textup{PI}}(x_k)\colon 0 \leqslant x_k < n_k\}$ in $O(n_k)$ time, yet provided $\lambda_k^*$ and $\phi_k^*$ are available. 

We next give a new result showing that, for the case in Proposition \ref{pro:osallocm3} of equal buffer  ($n_k \equiv n$) and server-pool ($m_k \equiv m$)
sizes at all queues, 
the PI index policy is insensitive to service rates $\mu_k$, reducing to SQ routing.

\begin{proposition}
\label{pro:model2simpl}
If $n_k \equiv n$ and $m_k \equiv m$ for each $k$, then the index
$\theta_k^{\textup{PI}}(x_k)$ is equal for each $k$, and hence the PI
 policy reduces to the SQ routing rule.
\end{proposition}
\begin{proof}
Proposition \ref{pro:osallocm3} ensures that, for each queue $k$, 
 $\lambda_k^* = 
\lambda \mu_k / \sum_l \mu_l$, and hence 
$r_k^* \equiv r^* \triangleq \lambda/ \sum_l \mu_l$. Now, the formulae
for the PI index in (\ref{eq:mipii1})--(\ref{eq:mipii2b}) ensure
that $\theta_k^{\textup{PI}}(x_k)$ is the same for each $k$, which yields the result. 
\end{proof}

\section{Numerical Experiments}
\label{s:bs}
This section reports on the results of a  small scale numerical study to
benchmark the proposed  RB second-order index policy, both against the optimal performance and against four alternative routing index policies that have been proposed in the literature: SQ, SED, PI, and the \emph{Never Queue} (NQ) policy discussed in \citet{shenkerWein89}.
Note that the routing indices corresponding to the SQ and the SED policies are given, respectively, by $\theta_k^{\textup{SQ}}(x_k) \triangleq x_k$, and 
$\theta_k^{\textup{SED}}(x_k) \triangleq 1/\mu_k$ if $0 \leqslant x_k < m_k$, with $\theta_k^{\textup{SED}}(x_k) \triangleq 1/\mu_k + (x_k + 1 - m_k)/(m_k \mu_k) = (x_k + 1)/(m_k \mu_k)$ if $m_k \leqslant x_k  < n_k$. Recall that
$x_k$ is the number of jobs present in queue $k$,  including the ones getting
service.

As for the NQ policy, which is a myopic policy that seeks to maximize the instantaneous throughput rate, it is defined in \citet{shenkerWein89} for the case of single-server queues as follows: ``The NQ policy chooses the fastest server that has
an empty queue; if there are no empty queues, the queue with
minimal $x_k/\mu_k$ is selected.''
Noting that the term $x_k/\mu_k$ is the expected remaining work and, also, the expected wait prior to entering service under FCFS for a job joining queue $k$ when it holds $x_k$ jobs, waiting or being served, we can generalize such an NQ policy to the case of a queue with $m_k$ servers, as follows: choose
the queue with the fastest servers that has some idle server; if there are no idle servers at any queue, select the queue with the
minimal expected wait prior to entering service under FCFS, which is given by $(x_k-m_k+1)/(m_k \mu_k)$.
The NQ index for queue $k$ can hence be formulated as $\theta_k^{\textup{NQ}}(x_k) \triangleq 1/\mu_k$ if $0 \leqslant x_k < m_k$, with $\theta_k^{\textup{NQ}}(x_k) \triangleq  c + (x_k + 1 - m_k)/(m_k \mu_k)$ if $m_k \leqslant x_k  < n_k$, where $c \triangleq \max_k 1/\mu_k$.

The study further assesses  two lower bounds on the minimum loss performance: that resulting from the relaxation  in the 
RB index method, which we denote by LBR,  and the  lower bound obtained by \emph{pooling} all buffer and service resources into a single-buffer single-server queue, which we denote by LBP.

The study considers instances with $K = 3$ queues.
For each base instance, the nominal load $\rho \triangleq \lambda/\sum_k m_k \mu_k$ was made to range over the interval $[0.7, 1.2]$,  from moderate load to overload, by varying the arrival rate $\lambda$.
For each  instance, the minimum loss probability ($z^{\textup{OP}}$) was computed by solving with CPLEX the LP formulation of the optimality equations for the problem obtained from (\ref{eq:pa2a}) by dropping the factor $\lambda$ from the objective. Such equations are: for every system state $\mathbf{0} \leqslant \mathbf{x} \leqslant \mathbf{n}$, 
\begin{equation}
\label{eq:zopdpeq}
z^{\textup{OP}} + h(\mathbf{x}) = 
\begin{cases}
\displaystyle \min_{k \in \mathbb{K}\colon x_k < n_k} \, -\sum_{l \in \mathbb{K}\colon 1 \leq x_l \leqslant n_k} \bar{\mu}_l(x_l) \Delta_l h(\mathbf{x}) + \lambda \Delta_k h(\mathbf{x}+\mathbf{e}_k), & \text{if } \mathbf{x} \neq \mathbf{n} \\ 
\displaystyle 1 -\sum_{l \in \mathbb{K}\colon 1 \leq x_l \leqslant n_k} \bar{\mu}_l(x_l) \Delta_l h(\mathbf{x}), & \text{if } \mathbf{x} = \mathbf{n},
\end{cases}
\end{equation}
where the $h(\mathbf{x})$ represent  the \emph{relative costs}, and we write  $\Delta_l h(\mathbf{x}) \triangleq h(\mathbf{x}) - h(\mathbf{x}-\mathbf{e}_l)$, with $\mathbf{e}_l$ being the $l$th unit coordinate vector in $\mathbb{R}^K$.
 Further, the loss probabilities $z^{\textup{SQ}}$, $z^{\textup{SED}}$, $z^{\textup{PI}}$, $z^{\textup{NQ}}$, and $z^{\textup{RB}}$ under the SQ, SED, OBS, PI, NQ, and RB policies, were computed by solving with MATLAB the  corresponding evaluation equations.
Finally, the lower bounds $z^{\textup{LBP}}$ and $z^{\textup{LBR}}$ on $z^{\textup{OP}}$ were computed.

The results show both (a) the percent deviation from the minimum loss performance for each policy and bound considered, e.g., 
$100 (z^{\textup{PI}} - z^{\textup{OP}}) /  z^{\textup{OP}}$ and $100 (z^{\textup{OP}}-z^{\textup{LBP}}) /  z^{\textup{OP}}$; and 
(b) the percent performance improvement of the RB index policy against each of the other policies, e.g., $100 (z^{\textup{PI}} - z^{\textup{RB}}) /  z^{\textup{PI}}$.

The experiments focus on two highly asymmetric instances, as the author has found in other experiments  that it is on such instances that the largest differences between the performances of alternative policies are obtained.

\subsection{Experiment 1}
\label{s:exp1}
The base instance for the first experiment has server pool sizes $\mathbf{m} = (1, 4, 10)$, service rates $\boldsymbol{\mu} = (80, 15, 5)$, and buffer sizes $\mathbf{n} = (16, 12, 10)$.
Such a  buffer size dimensioning ensures that the largest expected delay (including the service time) of a job routed to any queue does not exceed $n_k/(m_k \mu_k) \equiv 0.2$ time units.

\begin{figure}[!htb]
\centering
\includegraphics[height=2.8in]{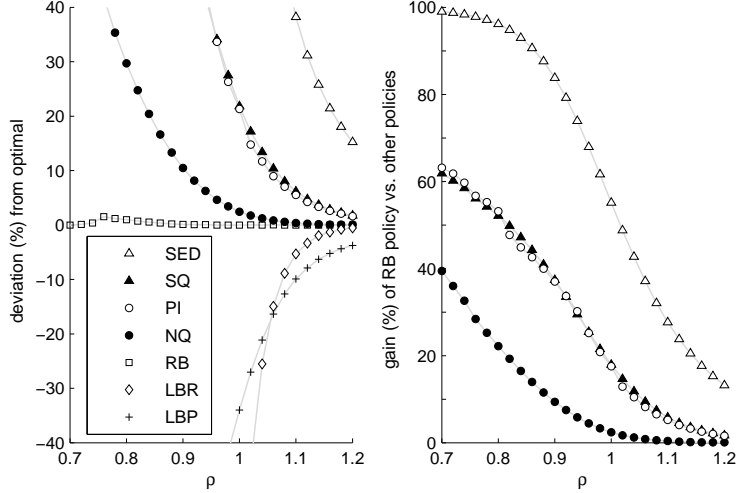}
\caption{Experiment 1: $\mathbf{m} = (1, 4, 10)$, $\boldsymbol{\mu} = (80, 15, 5)$, $\mathbf{n} = (16, 12, 10)$.}
\label{fig:exp1Nov11Plot}
\end{figure}

 Figure \ref{fig:exp1Nov11Plot} plots the results. 
The left pane shows that only the RB policy is near optimal, with the remaining policies being severely suboptimal in moderate traffic, being ranked in the order (from better to worse) NQ, PI and SQ (with PI and SQ having a similar performance), and SED.
The bounds are non-informative except under overload  ($\rho > 1$), with the LBR bound being better than LBP for large $\rho$.
 
The right pane shows that the gains of the RB policy against the alternative policies are substantial in moderate traffic, getting smaller as traffic becomes heavier.
The close performance of all policies under overload is to be expected, since in such conditions the system will be full most of the time, leaving little margin to the system controller for affecting performance.

\subsection{Experiment 2}
\label{s:exp2}
\begin{figure}[!htb]
\centering
\includegraphics[height=2.8in]{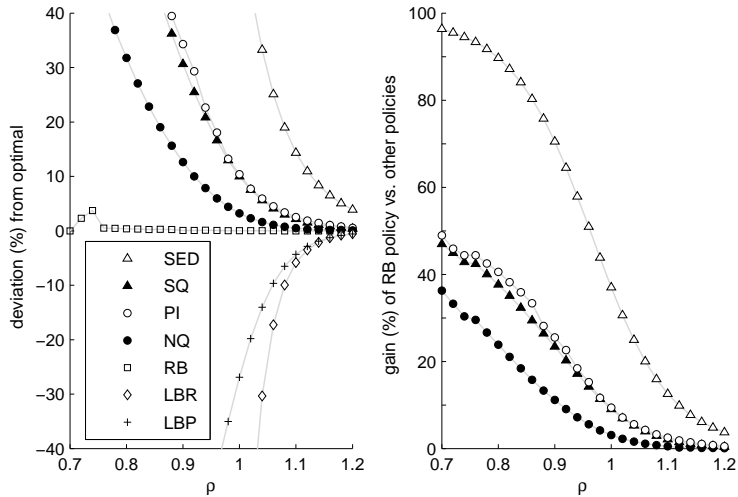}
\caption{Experiment 2: $\mathbf{m} = (1, 6, 8)$, $\boldsymbol{\mu} = (1440, 160, 100)$, $\mathbf{n} = (18, 12, 10)$.}
\label{fig:exp2Nov11Plot}
\end{figure}

The second experiment considers the base instance with server pool sizes $\mathbf{m} = (m_k) = (1, 6, 8)$, service rates $\boldsymbol{\mu} = (1440, 160, 100)$, and buffer sizes $\mathbf{n} = (18, 12, 10)$.
Such a  buffer size dimensioning ensures that the largest expected delay  of a job routed to any queue does not exceed $n_k/(m_k \mu_k) \equiv 0.0125$ time units.
 Figure \ref{fig:exp2Nov11Plot} plots the results, which are similar to those of the first experiment, except that now the bound LBP is always better than LBR. 

\subsection{Experiment 3}
\label{s:exp3}
\begin{figure}[!htb]
\centering
\includegraphics[height=2.8in]{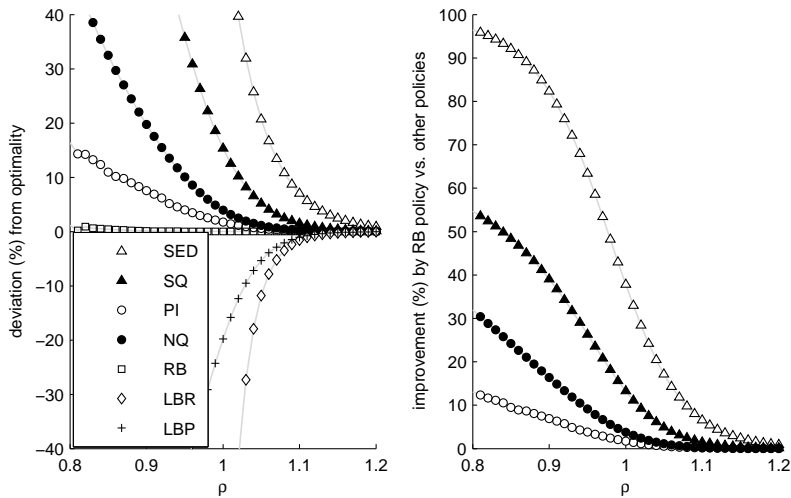}
\caption{Experiment 3: $\mathbf{m} = (1, 4, 10)$, $\boldsymbol{\mu} = (80, 15, 5)$, $\mathbf{n} = (18, 18, 18)$.}
\label{fig:exp3Feb12c}
\end{figure}

In experiments 1 and 2, model parameters were chosen to ensure that the largest expected delay of a job did not exceed a small value, motivated by QoS considerations. To explore the effect of allowing larger maximum expected delays, we report on a third experiment, which considers the base instance with server pool sizes and service rates as in the first experiment, but with the larger buffer sizes $\mathbf{n} = (18, 18, 18)$, under which the largest
 expected delay  of a job ($0.36$ time units) increases by $80\%$ with respect to that in the first experiment.
 Figure \ref{fig:exp3Feb12c} plots the results, which are similar to those of the first experiment, except that now the PI policy shows a substantially improved performance, becoming the second-best policy. Yet, the RB policy is still the best policy, being again nearly optimal and substantilly outperforming the other policies while the traffic intensity is not too high.

A similar behavior is observed when buffer sizes are modified in the second experiment to be $\mathbf{n} = (18, 18, 18)$.

Based on this and other experiments not reported here, the effect of increasing the buffer sizes on the performance of the tested policies (which has the effect of increasing the largest expected delay of a job) appears to be most apparent in the case of the PI policy, whose performance improves significantly. Yet, in all the instances considered by the author, the best policy among those considered is still the RB policy, which is consistently near optimal.

\section{Conclusions}
\label{s:concl}
This paper has presented a novel index policy for dynamic routing to parallel multiserver finite-buffer queues, with a maximum average throughput (or, equivalently, minimum average loss) objective, based on a  second-order RB index that is used to break ties caused by the fact that the conventional first-order RB index is constant. The proposed RB index is efficiently computed by a first-order linear recursion, and is also given in closed form. As a state-of-the-art alternative, the paper has considered the previously proposed and more computationally demanding index policy obtained by performing one step of the PI algorithm on a base policy given by the OBS. New results of algorithmic and theoretical interest are presented, including the result that the PI index policy reduces to SQ routing in the case of equal buffer and server-pool sizes.

A small scale numerical study is reported showing that, on the  instances considered, the proposed RB index policy is consistently near optimal and can substantially outperform each of four alternative index policies, which are the SQ, SED, NQ, and PI. 
The study, though small in scale, provides evidence for arguing that, in highly heterogeneous systems, it may not suffice with naive index policies to attain a near optimal throughput performance. 
The RB index policy, though relatively sophisticated, can nevertheless be efficiently evaluated in linear time on the number of index values to be computed, and is also given in closed form. 

It would be interesting for future work to try to establish theoretically the observed near-optimality of the RB index policy, or to identify a parameter range where it performs poorly, if any.
In the experiments performed by the author, including those not reported here, the RB index policy was found to be consistently near optimal.

\section*{Acknowledgements}
This work was partially supported by the Spanish Ministry of Education and Science project MTM2007-63140 and by the Spanish Ministry of Science and Innovation project MTM2010-20808.

\bibliographystyle{elsarticle-harv}



\end{document}